\documentclass[12pt]{amsart}
\usepackage{a4}

\DeclareMathAlphabet{\eusm}{U}{}{}{}  
\SetMathAlphabet\eusm{normal}{U}{eus}{m}{n}
\SetMathAlphabet\eusm{bold}{U}{eus}{b}{n}

\DeclareMathAlphabet{\eufrak}{U}{}{}{}  
\SetMathAlphabet\eufrak{normal}{U}{euf}{m}{n}
\SetMathAlphabet\eufrak{bold}{U}{euf}{b}{n}

\swapnumbers

\newtheorem{theorem}{Theorem}[section]

\newtheorem{corollary}[theorem]{Corollary}

\theoremstyle{definition}
\newtheorem{definition}[theorem]{Definition}
\newtheorem{example}[theorem]{Example}

\theoremstyle{remark}
\newtheorem{remark}[theorem]{Remark}
\numberwithin{equation}{section}

\begin{document}
\title{Free L\'evy Processes on Dual Groups}
\author{Uwe Franz}
\address{Institut f\"ur Mathematik und Informatik,
  Ernst-Moritz-Arndt-Universit\"at Greifs\-wald, Friedrich-Ludwig-Jahnstr. 15A,
  D-17487 Greifs\-wald, Germany,
\newline Email: {\tt franz@uni-greifswald.de},
\newline Homepage: {\tt
  http://hyperwave.math-inf.uni-greifswald.de/algebra/franz}}

\thanks{Extended abstract to be published in Mini-proceedings: Second
  MaPhySto Conference on ``L\'evy Processes - Theory and Applications,''
  January 2002, Aarhus, Denmark}

\begin{abstract}
We give a short introduction to the theory of L\'evy processes on dual
groups. As examples we consider L\'evy processes with additive increments and
L\'evy processes on the dual affine group.
\end{abstract}
\maketitle

\section{Introduction}

L\'evy processes play a fundamental r\^ole in probability theory and have many
important applications in other areas such as statistics, financial mathematics,
functional analysis or mathematical physics, as well.

In quantum probability they first appeared in a model for the laser in
\cite{waldenfels84}. This lead to the theory of L\'evy processes on involutive
bialgebras, cf.\ \cite{accardi+schuermann+waldenfels88,schuermann93,franz+schott99}. The
increments of these L\'evy processes are independent in the sense of tensor
independence, which is a straightforward generalisation of the notion of
independence used in classical probability theory. However, in quantum probability
there exist also other notions of independence like, e.g., freeness \cite{voiculescu+dykema+nica92}, see Paragraph
\ref{par-ind}. In order to formulate a general theory of L\'evy processes for
these independences, the $*$-bialgebras or quantum groups have to be replaced by
the dual groups introduced in \cite{voiculescu87}, see
\cite{schuermann95b,benghorbal+schuermann99,franz01,franz01b}.

In this paper we give an introduction to the theory of L\'evy processes
on dual groups, which avoids most of the algebraic prerequisites. In
particular, we will not define dual groups, but only consider two examples,
namely tensor and free L\'evy processes with additive increments and tensor
and free L\'evy processes on the
dual affine group. Our approach is similar to rewriting the definition of
classical Lie group-valued L\'evy processes in terms of a coordinate system, see
Definitions \ref{def-add-ten-levy}, \ref{def-add-free-levy},
\ref{def-aff-ten-levy}, and \ref{def-aff-free-levy}.

Quantum L\'evy processes play an important r\^ole in the theory of continuous
measurement, cf. \cite{holevo01}, and in the theory of dilations, where they
describe the evolution of a big system or heat bath, which is coupled to the
small system whose evolution one wants to describe.

Additive free L\'evy processes where first studied in
\cite{glockner+schuermann+speicher92}, and more recently in \cite{biane98,anshelevich01b,anshelevich01c,barndorff-nielsen+thorbjornsen01a,barndorff-nielsen+thorbjornsen01b}.

\section{Preliminaries}

In this section we introduce the basic notions and definitions that we will use. For more detailed introductions to quantum
probability see, e.g., \cite{parthasarathy92,biane93,meyer95,holevo01}.

\subsection{Quantum probability}

Non-commutative probability or quantum probability is motivated by the
statistical interpretation of quantum mechanics where an operator is
interpreted as an analog of a random variable. The r\^ole of the classical
probability space is played by a (pre-)Hilbert space $\mathcal{H}$ and the
measure is replaced by a unit vector $\Omega\in \mathcal{H}$ called state vector.

In this paper we will mean by a (real) quantom random variable $X$ on
$(\mathcal{H},\Omega)$ a (symmetric)
linear operator on the pre-Hilbert space $\mathcal{H}$, which has an adjoint,
i.e.\ for which there exists a linear operator $X^*$, such that
\[
\langle u, Xv\rangle = \langle X^*u,v\rangle
\]
for all $u,v\in\mathcal{H}$. Its law (w.r.t.\ the state
vector $\Omega$) is the functional
$\phi_X:\mathbb{C}[x]\to\mathbb{C}$ on the algebra $\mathbb{C}[x]$ of
polynomials in one variable defined by
\[
\phi_X(x^k)=\langle\Omega,X^k\Omega\rangle,
\]
for $k\in\mathbb{N}$. If $X$ is symmetric, then there exists a (possibly
non-unique) probability measure $\mu$ on $\mathbb{R}$ such that
\[
\phi_X(x^k)=\int_\mathbb{R}x^k{\rm d}\mu.
\]

Let $X$ be a classical $\mathbb{R}$- or $\mathbb{C}$-valued random variable
with finite moments on some probability space $(M,\mathcal{M},P)$. It
becomes a quantum random variable on
$\mathcal{H}=L^\infty(M,\mathcal{M},P)$, if we let it act on the bounded
functions on $M$ 
by multiplication, $L^\infty(M)\ni f\mapsto Xf\in L^\infty(M)$ with
$Xf(m)=X(m)f(m)$ for $m\in M$. If we take the constant function $\Omega(m)=1$
for all $m\in M$ for the state vector, then we recover the classical
distribution of $X$, i.e.,
\[
\langle\Omega, X^k\Omega\rangle = \int_M X^k{\rm d}P=\mathbb{E}(X^k),
\]
for $k\in\mathbb{N}$. If $X$ is $\mathbb{R}$-valued, then it is also real as
quantum random variable, i.e.\ symmetric.

A (real) quantum random vector on $(\mathcal{H},\Omega)$ is an $n$-tuple
$X=(X_1,\ldots,X_n)$ of (real) quantum random
variables on $(\mathcal{H},\Omega)$. Its law is the functional $\phi_X:\mathbb{C}\langle
x_1,\ldots,x_n\rangle\to\mathbb{C}$ on the algebra of non-commutative
polynomials $\mathbb{C}\langle x_1,\ldots,x_n\rangle$ in $n$ variables defined
by
\[
\phi_{X}( x_{i_1}^{k_1}\cdots x_{i_r}^{k_r}) = \langle\Omega,
X_{i_1}^{k_1}\cdots X_{i_r}^{k_r}\Omega\rangle,
\]
for all $i_1,\ldots,i_r\in\{1,\ldots,n\}$, $k_1,\ldots,k_r\in\mathbb{N}$.

A (real) operator process is an indexed family $(X_\imath)_{\imath\in I}$ on
$(\mathcal{H},\Omega)$ of (real) quantum
random variables or (real) quantum random vectors on
$(\mathcal{H},\Omega)$. We will call two operator processes
$(X_\imath)_{\imath\in I}$ and $(Y_\imath)_{\imath\in I}$ equivalent, if they
have the same joint moments, i.e. if
\[
\langle\Omega,
X^{k_1}_{\imath_1}\cdots X^{k_r}_{\imath_r}\Omega\rangle=\langle\Omega,
Y^{k_1}_{\imath_1}\cdots Y^{k_r}_{\imath_r}\Omega\rangle
\]
for all $\imath_1,\ldots,\imath_r\in I$ and all $k_1,\ldots,k_r\in\mathbb{N}$.

\subsection{Freeness and Independence}\label{par-ind}

Let now $\mathcal{A}_1,\ldots,\mathcal{A}_k$ be algebras of adjoin\-table linear
operators on some pre-Hilbert space $\mathcal{H}$, closed under taking
adjoints and containing the identity operator $\mathbf{1}$.
\begin{definition}
$\mathcal{A}_1,\ldots,\mathcal{A}_k$ are called tensor independent (w.r.t.\ to
the state vector $\Omega$), if
\begin{itemize}
\item[(i)]
for all $1\le i,j\le k$ and all $X\in\mathcal{A}_i$ and $Y\in\mathcal{A}_j$,
we have
\[
[X,Y]:=XY-YX=0,
\]
\item[(ii)]
and for all $X_1\in\mathcal{A}_1,\ldots,X_k\in\mathcal{A}_k$ we have
\[
\langle\Omega,X_1\cdots X_k\Omega\rangle = \langle\Omega,X_1\Omega\rangle
\cdots \langle\Omega,X_k\Omega\rangle.
\]
\end{itemize}
\end{definition}

This definition is the natural analogue of the notion of independence in
classical probability to our setting. It is also the one used in quantum physics when
one speaks of independent observables. But in quantum probability there exist
other, inequivalent notions of independence.

\begin{definition}
$\mathcal{A}_1,\ldots,\mathcal{A}_k$ are called free, if for all $1\le
i_1,\ldots,i_r\le k$ with $i_1\not= i_2\not=\cdots\not= i_r$ (i.e.\
neighboring indices are different) and all
$X_1\in\mathcal{A}_{i_1},\ldots,X_r\in\mathcal{A}_{i_r}$ with
\[
\langle\Omega,X_1\Omega\rangle = \cdots = \langle\Omega,X_r\Omega\rangle = 0,
\]
we have
\[
\langle\Omega,X_1\cdots X_r\Omega\rangle=0.
\]
\end{definition}

Quantum random variables or quantum random vectors $X,Y,Z,\ldots$ are called
tensor independent or free, iff the unital $*$-algebras they generate are tensor
independent or free.

\begin{remark}
These definitions allow to compute arbitrary joint moments of tensor independent
or free random variables from their marginal distributions.

For the free case this computation can be done recursively on the order of the
moment be expanding
\[
0= \varphi\Big(\big(X_1-\varphi(
X_1)\mathbf{1}\big)\cdots\big(X_r-\varphi(X_r)\mathbf{1}\big)\Big),
\]
where we wrote $\varphi(\cdot)$ instead of $\langle\Omega,\cdot\,\Omega\rangle$
for the expectation.

Let $X_1$ and $X_2$ be two free quantum random variables, then one obtains in this way
\[
0= \varphi\Big(\big(X_1-\varphi(
X_1)\mathbf{1}\big)\big(X_2-\varphi(X_2)\mathbf{1}\big)\Big) =
\varphi(X_1X_2)- \varphi(X_1)\varphi(X_2),
\]
and therefore
\[
\varphi(X_1X_2)=\varphi(X_1)\varphi(X_2).
\]
as for tensor independent quantum random variables or independent classical
random variables. But for higher moments the
formulas are different, one gets, e.g.,
\[
\varphi(X_1X_2X_1X_2)=\varphi(X_1^2)\varphi(X_2)^2+\varphi(X_1)^2\varphi(X_2^2)-\varphi(X_1)^2\varphi(X_2)^2.
\]
This formula can also be used to show that there exist no non-trivial examples
of commuting free quantum random variables. If $X_1$ and $X_2$ commute, then
we would get
\[
\varphi(X_1X_2X_1X_2)=\varphi(X_1^2X_2^2)=\varphi(X_1^2)\varphi(X_2)^2,
\]
since $X_1^2$ and $X_2^2$ are also free. Therefore
\[
\varphi\big(X_1-\varphi(X_1)\mathbf{1}\big)\varphi\big(X_2-\varphi(X_2)\mathbf{1}\big)=0,
\]
i.e.\ at least one of the two quantum random variables has a trivial distribution.
\end{remark}

\begin{remark}
Besides tensor independence and freeness there exist other notions of
independence that are used in quantum probability. In a series of papers
\cite{schuermann95,speicher96,benghorbal+schuermann99,muraki01b,muraki02} it
was shown that there exist exactly five ``universal'' notions of independence
satisfying a natural set of axioms.  Besides tensor independence and freeness
these are boolean, monotone, and anti-monotone independence. In \cite{franz01b}
the boolean, monotone, and anti-monotone independence where reduced to tensor
independence. If this is also possible for freeness is still an open problem.
\end{remark}

\section{Additive L\'evy Processes}

\begin{definition}\label{def-add-ten-levy}
An operator process $(X_t)_{t\ge0}$ on $(\mathcal{H},\Omega)$ is called an
additive tensor L\'evy process (w.r.t.\ $\Omega$), if the increments
\[
X_{st}:=X_t-X_s,
\]
are
\begin{itemize}
\item[(i)]
tensor independent, i.e.\ the quantum random variables $X_{s_1t_1},\ldots,X_{s_rt_r}$ are tensor
independent for all $0\le s_1\le t_1\le s_2\le \cdots\le s_r\le t_t$,
\item[(ii)]
stationary, i.e.\ the law of an increment depends only on $t-s$, and
\item[(iii)]
weakly continuous, i.e.
$\lim_{t\searrow s} \langle\Omega,X_{st}^k\Omega\rangle=0$
for $k=1,2,\ldots$.
\end{itemize}
\end{definition}

Replacing tensor independence by another universal notion of independence we
can define the corresponding other classes of L\'evy processes. E.g., for
freeness we get the following definition.

\begin{definition}\label{def-add-free-levy}
An operator process $(X_t)_{t\ge0}$ on $(\mathcal{H},\Omega)$ is called an
additive free L\'evy process (w.r.t.\ $\Omega$), if the increments
\[
X_{st}:=X_t-X_s,
\]
are
\begin{itemize}
\item[(i')]
free, i.e.\ the quantum random variables
$X_{s_1t_1},\ldots,X_{s_rt_r}$ are free for all $0\le s_1\le t_1\le s_2\le \cdots\le s_r\le t_t$,
\end{itemize}
and satisfy conditions (ii) and (iii) of Definition \ref{def-add-ten-levy}.
\end{definition}

For each of these notions of independence one can define a Fock space and
creation, annihilation, and conservation or gauge operators on this Fock space.

For example the (algebraic) free Fock space
$\mathcal{H}=\mathcal{F}_{\rm F}(\eufrak{h})$ over a (pre-)Hilbert space $\eufrak{h}$ is defined as
\[
\mathcal{F}_{\rm F}(\eufrak{h})= \bigoplus_{n=0}^\infty \eufrak{h}^{\otimes n}
\]
where $\eufrak{h}^{\otimes 0}=\mathbb{C}$. The vector $\Omega=1+0+\cdots$ is
called the vacuum vector. For  a vector $u\in\eufrak{h}$ we can define the
creation operator $a^+(u)$ and the annihilation operator $a^-(u)$ by
\begin{eqnarray*}
a^+(u) u_1\otimes\cdots u_k &=& u\otimes u_1\otimes \cdots \otimes u_k, \\
a^-(u) u_1\otimes\cdots u_k &=& \langle u, u_1\rangle u_2\otimes \cdots\otimes
u_k.
\end{eqnarray*}
These operators are mutually adjoint, on the vacuum vector they act as
$a^+(u)\Omega=u$ and $a^-(u)\Omega=0$.

The conservation operator $\Lambda(X)$ of some linear operator $X$ on
$\eufrak{h}$ is defined by
\[
\Lambda(X)u_1\otimes\cdots u_k = (Xu_1)\otimes u_2\otimes\cdots\otimes u_k
\]
and $\Lambda(X)\Omega=0$. It satisfies $\Lambda(X)^*=\Lambda(X^*)$.

Glockner, Sch\"urmann, and Speicher have shown that every additive free L\'evy
process can be realized as a linear combination of these three operators and time.

\begin{theorem}\label{theo-add-levy}
\cite{glockner+schuermann+speicher92}
Let $(X_t)_{t\ge 0}$ be an additive free L\'evy process. Then there exists a
pre-Hilbert space $\eufrak{k}$, a linear operator $T$ on $\eufrak{k}$, vectors
$u,v\in\eufrak{k}$, and a scalar $\lambda\in\mathbb{C}$ such that $(X_t)_{t\ge
  0}$ is equivalent to the operator process $(X'_t)_{t\ge 0}$ on the free Fock
space $\mathcal{F}_{\rm F}\big(L^2(\mathbb{R}_+,\eufrak{k})\big)$ over
$\eufrak{h}=L^2(\mathbb{R}_+,\eufrak{k})\cong L^2(\mathbb{R}_+)\otimes\eufrak{k}$ defined by
\[
X'_t = \Lambda(\chi_{[0,t]}\otimes T) + a^+(\chi_{[0,t]}\otimes u) +
a^-(\chi_{[0,t]}\otimes v) + t\lambda\mathbf{1}
\]
for $t\ge 0$. Furthermore, if we require that $\eufrak{k}$ is spanned by
$\{T^ku,T^kv|k=0,1,\ldots\}$, then $\eufrak{k}$, $T$, $u$, $v$, and $\lambda$ are
unique up to unitary equivalence.

$(X'_t)_{t\ge 0}$ is symmetric, if and only if $T^*=T$, $u=v$ and
$\lambda\in\mathbb{R}$ in the unique minimal tuple.
\end{theorem}

\begin{remark}
Analogous results hold for the other universal independences. For tensor
independence see \cite{schuermann91c}, for the boolean \cite{benghorbal01}, and
\cite{franz01b} for the monotone case. Note that in the boolean and in the
monotone case the time process has to be modified.
\end{remark}

The five independences can also be used to define convolutions for compactly
supported measures. Let $\mu_1$ and $\mu_2$ be two compactly supported
probability measures on $\mathbb{R}$ and choose two independent real quantum random
variables $X_1$ and $X_2$ on some pre-Hilbert space $\mathcal{H}$ such that
\[
\langle\Omega, X^k_i\Omega\rangle = \int_\mathbb{R} x^k{\rm d}\mu_i
\]
for all $k\in \mathbb{N}$ and $i=1,2$ and some unit vector
$\Omega\in\mathcal{H}$. The operator $X_1+X_2$ is again symmetric and
bounded, therefore there exists a unique compactly supported probability
measure $\mu$ such that
\[
\int_\mathbb{R}x^k{\rm d}\mu = \langle\Omega, (X_1+X_2)^k\Omega\rangle
\]
for all $k\in \mathbb{N}$

It is always possible to construct such a pair and the law of $X_1+X_2$
depends only on the laws of $X_1$ and $X_2$ and the notion of independence that
has been chosen.

If $X_1$ and $X_2$ are tensor independent, then the measure $\mu$ obtained
in this way is the usual additive convolution of $\mu_1$ and $\mu_2$. If $X_1$
and $X_2$ are free, then $\mu$ is the free additive convolution of  $\mu_1$
and $\mu_2$.

These convolutions can actually be defined for arbitrary  probability
measures. 

It is possible to show that infinitely divisible measures can be embedded into
continuous convolution semigroups in all five cases and that furthermore there
exists a L\'evy process for every continuous convolution semigroup. This shows
that in all five cases the infinitely divisible measures on $\mathbb{R}$
(which are characterized by their moments) can be classified by tuples $(\eufrak{k},T,u,\lambda)$ consisting of a pre-Hilbert
space $\eufrak{k}$, a symmetric operator $T$ on $\eufrak{k}$, a vector
$u\in\eufrak{k}$, and a real number $\lambda$.
\begin{corollary}
There exist bijections (up to moment uniqueness) between the five classes of
infinitely divisible measures with finite moments.
\end{corollary}
\begin{remark}
The bijection between the usual infinitely divisible measures and the freely
infinitely divisible measures is known under the name Pata-Bercovici
bijection, cf.\ \cite{bercovici+pata99}, it actually extends to all infinitely
divisible measures, not just those characterized by their moments, and has
many useful properties, cf.\ \cite{barndorff-nielsen+thorbjornsen01a} and the
references therein.

For example, the Bercovici-Pata bijection is a homomorphism between the usual
infinitely divisible measures and the freely infinitely divisible measures and
their respective convolutions. This is not the case for the bijection between  usual
infinitely divisible measures and the monotone infinitely divisible measures,
because due to the non-commutativity of the monotone convolution this is
impossible. For the L\'evy-Khintchine formula for the boolean and monotone
case, see \cite{speicher+woroudi93,muraki00}.
\end{remark}

\begin{definition}\label{def-gauss-poisson}
Let $(X_t)_{t\ge 0}$ be a real additive L\'evy process for one of the five
universal independences.

If there exists a tuple $(\eufrak{k},T,u,\lambda)$ for $(X_t)_{t\ge 0}$ with $T=0$, then $(X_t)_{t\ge 0}$ is called Gaussian.

If there exists a tuple $(\eufrak{k},T,u,\lambda)$ for $(X_t)_{t\ge 0}$ and
vector $\omega\in\eufrak{k}$ such that $u=T\omega$ and
$\lambda=\langle\omega,T\omega\rangle$, then $(X_t)_{t\ge 0}$ is called a
compound Poisson process. 
\end{definition}

If $(X_t)_{t\ge 0}$ is Gaussian, then the unique minimal tuple associated to it by
Theorem \ref{theo-add-levy} has the form $(\mathbb{C}, 0, z,\lambda)$ and
$(X_t)_{t\ge 0}$ can be realized as a sum of creation, annihilation and time
only, with no conservation part.

\begin{example}
Let $(X_t)_{t\ge 0}$ be a classical compound Poisson process with L\'evy
measure $\mu$, i.e. with characteric function
\[
\mathbb{E}\left(e^{iuX_t}\right) = \exp\left(t\int_{\mathbb{R}\backslash\{0\}}
  (e^{iux}-1){\rm d}\mu(x)\right).
\]
We assume that $\mu$ has finite moments, then $(X_t)_{t\ge 0}$ is an additive
tensor L\'evy process in the sense of Definition \ref{def-add-ten-levy}.
We can define a tuple $(\eufrak{k},T,u,\lambda)$ for
$(X_t)_{t\ge 0}$ by Theorem \ref{theo-add-levy} as follows. For the
pre-Hilbert space $\eufrak{k}$ we take the space of polynomials
\[
\eufrak{k}={\rm span}\,\{x^k;k=0,1,2,\ldots\},
\]
considered as a subspace of the Hilbert space $L^2(\mathbb{R},\mu)$, i.e.,
with the inner product
\[
\langle x^k,x^\ell\rangle = \int_{\mathbb{R}} x^{k+\ell}{\rm d}\mu(x),
\]
and divided by the the nullspace of this inner product, if $\mu$ is finitely
supported. The operator $T$ is
multiplication by $x$, i.e., $Tx^k = x^{k+1}$, the vector $u$ is the function
$f(x)=x$, and the scalar $\lambda$ is the first moment of $\mu$, i.e.,
$\lambda=\int_{\mathbb{R}} x{\rm d}\mu(x)$.

Taking for $\omega$ the constant function $1$, we see that $(X_t)_{t>0}$ is
also a compound Poisson process in sense of Definition \ref{def-gauss-poisson}.
\end{example}

Theorem \ref{theo-add-levy} can also be used to give an
It\^o-L\'evy-type decomposition of additive quantum L\'evy processes.

\begin{corollary}
Let $(X_t)_{t\ge 0}$ be a real additive L\'evy process for one of the five
universal independences. Then $(X_t)_{t\ge 0}$ can be realized as a sum of a
Gaussian L\'evy process $(X^{\rm G}_t)_{t\ge 0}$ and a ``jump'' part $(X^{\rm P}_t)_{t\ge 0}$, which
can be approximated by (compensated) compound Poisson processes.
\end{corollary}
\begin{proof}
We only briefly outline the proof.

Let $(\eufrak{k},T,u,\lambda)$ be a tuple for $(X_t)_{t\ge 0}$. Since $T$ is
symmetric, we can decompose the closure of $\eufrak{k}$ into a direct sum of
the closures of the kernel of $T$ and the image of $T$. Let $u=u_0+u_1$ with
$u_0\in\overline{{\rm ker}\,T}$ and $u_1\in\overline{{\rm im}\,T}$.

If $u_1$ is actually in the image of $T$, then there exists a vector
$\omega\in\eufrak{k}$ with $T\omega=u_1$ and $(X_t)_{t\ge 0}$ is equivalent to the sum
of the Gaussian L\'evy process $(X^{\rm G}_t)_{t\ge 0}$ with the tuple
$(\mathbb{C},0,u_0,\lambda-\langle \omega,T \omega\rangle)$ and the
compound Poisson process $(X^{\rm P}_t)_{t\ge 0}$ with the tuple $(\eufrak{k},T,u_1,\langle \omega,T\omega\rangle)$.

If $u_1$ is not in the image of $T$, then we take a sequence
$\omega_n\in\eufrak{k}$ such that $\lim T\omega_n=u_1$ and define the ``jump'' part by
\begin{eqnarray*}
X^{\rm P}_t &=&\Lambda(\chi_{[0,t]}\otimes T)+ a^+(\chi_{[0,t]}\otimes u_1)+a^-(\chi_{[0,t]}\otimes u_1) \\
&=& \lim_{n\to\infty} \Big(\chi_{[0,t]}\otimes \Lambda(T)+ a^+\big(\chi_{[0,t]}\otimes (T\omega_n)\big)+a^-\big(\chi_{[0,t]}\otimes (T\omega_n)\big)\Big),
\end{eqnarray*}
i.e.\ as the limit of compensated compound Poisson processes. The Gaussian
part is then determined by the tuple $(\mathbb{C},0,u_0,\lambda)$.
\end{proof}
\begin{remark}
Using the spectral representation $\overline{T}=\int_\mathbb{R} x{\rm d}P_x$
of the closure of $T$,
the ``jump'' part can be written as an integral over the ``jump'' sizes.

However, note that the It\^o-L\'evy-type decomposition gives a decomposition into a
continuous Gaussian part and a jump part only in the tensor case. In the other
cases the classical processes that can be associated to the corresponding
Gaussian processes do not have continuous paths, see, e.g., \cite{biane98}.

Using different methods and not assuming the existence of moments,
Barndorff-Nielsen and Thorbj{\o}rnson \cite{barndorff-nielsen+thorbjornsen01b}
have also obtained an It\^o-L\'evy decomposition for additive free L\'evy processes.
\end{remark}

\section{L\'evy Processes on the (Dual) Affine Group}

Recall that the affine group can be defined as the group of matrices
\[
{\rm Aff}=\left\{\left(\begin{array}{cc}  a & b \\ 0 & 1\end{array}\right) : a >0,
    b\in\mathbb{R}\right\}.
\]
The calculation
\[
\left(\begin{array}{cc}  a_1 & b_1 \\ 0 & 1\end{array}\right)\left(\begin{array}{cc}  a_2
      & b_2 \\ 0 & 1\end{array}\right) = \left(\begin{array}{cc}  a_1a_2 & a_1b_2+b_1 \\ 0 &
      1\end{array} \right)
\]
shows that the group multiplication takes the form
\begin{eqnarray*}
A(g_1g_2) &=& A(g_1)A(g_2), \\
B(g_1g_2) &=& A(g_1)B(g_2)+B(g_1),
\end{eqnarray*}
for the coordinates $A,B$ defined by
\[
A\left(\begin{array}{cc}  a & b \\ 0 & 1\end{array}\right) = a, \qquad
B\left(\begin{array}{cc}  a & b \\ 0 & 1\end{array}\right) = b.
\]

We define tensor L\'evy processes on the dual affine group in term of
increments. Since the increments are tensor independent for disjoint time
intervals, they commute, and we can write the products in the multiplication
formula in any order we like.
\begin{definition}\label{def-aff-ten-levy}
An operator process $\big((A_{st},B_{st})\big)_{0\le s\le t}$ on
$(\mathcal{H},\Omega)$ is called a (left) tensor L\'evy process on the dual affine
group (w.r.t.\ $\Omega$), if the following four conditions are satisfied.
\begin{itemize}
\item[(i)]
(Increment property) For all $0\le s\le t\le u$,
\begin{eqnarray*}
A_{su} &=& A_{tu}A_{st}, \\
B_{su} &=& A_{tu}B_{st} + B_{tu}.
\end{eqnarray*}
\item[(ii)]
(Independence) The increments
$(A_{s_1t_1},B_{s_1t_1}),\ldots,(A_{s_rt_r},B_{s_rt_r})$ are tensor
independent for all $0\le s_1\le t_1\le s_2\le \cdots\le s_r\le t_r$.
\item[(iii)]
(Stationarity) The law of $(A_{st},B_{st})$ depends only on $t-s$.
\item[(iv)]
(Weak continuity) For all $k_1,\ldots,k_r,\ell_1,\ldots,l_r\in\mathbb{N}$, we
have
\[
\lim_{t\searrow s} \langle\Omega,A_{st}^{k_1} B_{st}^{\ell_1}\cdots
A_{st}^{k_r} B_{st}^{\ell_r}\Omega\rangle=
\left\{
\begin{array}{rcr}
1 & \mbox{ if }& \ell_1+\cdots \ell_r=0, \\
0 & \mbox{ if }& \ell_1+\cdots \ell_r>0.
\end{array}\right.
\]
\end{itemize}
\end{definition}
Every L\'evy process with values in the semi-group
\[
\left\{\left(\begin{array}{cc}  a & b \\ 0 & 1\end{array}\right) : a, b\in\mathbb{R}\right\}
\]
with finite moments gives an example of a tensor L\'evy process on the dual
affine group in the sense of our definition, since we didn't impose that the
$A_{st}$ are strictly positive or invertible.

\begin{example}
There are also examples in which $A_{st}$ and $B_{st}$ do not commute and which
do not correspond to classical L\'evy processes. E.g., the quantum Az\'ema
martingale \cite{emery89,parthasarathy90,schuermann91b} with paramater $q\in\mathbb{R}$ defined by the quantum
stochastic differential equations
\begin{eqnarray*}
{\rm d}A_{st} &=& A_{st}{\rm d}\Lambda_t(q-1), \\
{\rm d}B_{st} &=& B_{st}{\rm d}\Lambda_t(q-1) + {\rm d}a_t^+(1) + {\rm d}a_t^-(1),
\end{eqnarray*}
on the Bose Fock space $\mathcal{F}_{\rm T}(L^2(\mathbb{R}_+)$, with initial
 conditions
\begin{eqnarray*}
A_{ss} &=& {\rm id}, \\
B_{ss} &=& 0,
\end{eqnarray*}
defines a tensor L\'evy process on the dual affine group. For $q=1$, we have
$A_{st}={\rm id}$ for all $0\le s\le t$ and in the vacuum state $(B_{st})$ is
equivalent to classical Brownian motion. For $q\not=1$,
$A_{st}$ and $B_{st}$ do not commute and $(B_{st})$ is equivalent to the
classical Az\'ema martingale with parameter $c=q-1$.
\end{example}

When we want to define free L\'evy processes, different orders of the products
in the multiplication formula will lead to different classes of L\'evy
processes. The choice in the definition proposed here is motivated by the fact
that if $B_{st}$ and $B_{tu}$ are symmetric, then $B_{su}$ is also symmetric.
\begin{definition}\label{def-aff-free-levy}
An operator process $\big((a_{st},B_{st})\big)_{0\le s\le t}$ on
$(\mathcal{H},\Omega)$ is called a (left) free L\'evy process on the dual affine
group (w.r.t.\ $\Omega$), if the conditions
\begin{itemize}
\item[(i')]
(Increment property) For all $0\le s\le t\le u$,
\begin{eqnarray*}
a_{su} &=& a_{tu}a_{st}, \\
B_{su} &=& a_{tu}B_{st}a^*_{tu}+B_{st}.
\end{eqnarray*}
\item[(ii')]
(Independence) The increments
$(a_{s_1t_1},B_{s_1t_1}),\ldots,(a_{s_rt_r},B_{s_rt_r})$ are free for all $0\le s_1\le t_1\le s_2\le \cdots\le s_r\le t_r$.
\end{itemize}
and conditions (iii) and (iv) from the previous definition are satisfied (with
$A_{st}=a_{st}a^*_{st}$).
\end{definition}
Note that with this definition $A_{st}=a_{st}a^*_{st}$ is automatically positive.

\begin{example}
Let $\gamma\in\mathbb{C}$. The operator process $\big((a_{st},B_{st})\big)_{0\le s\le t}$ defined by the
quantum stochastic equations
\begin{eqnarray*}
{\rm d}a_{st} &=& {\rm d}\Lambda_t(\gamma-1)a_{st}, \\
{\rm d}B_{st} &=&  {\rm d}\Lambda_t(\gamma-1)B_{st} + B_{st}{\rm d}\Lambda_t(\overline{\gamma}-1) + {\rm d}a_t^+(1) + {\rm d}a_t^-(1),
\end{eqnarray*}
on the free Fock space  $\mathcal{F}_{\rm F}(L^2(\mathbb{R}_+)$, with initial
 conditions
\begin{eqnarray*}
a_{ss} &=& {\rm id}, \\
B_{ss} &=& 0,
\end{eqnarray*}
defines a free L\'evy process on the dual affine group. For $\gamma=1$, we get
$a_{st}={\rm id}$ and $(B_{st})$ is equal to the free Brownian motion,
\[
B_{st}=a^+(\chi_{[s,t[})+a^-(\chi_{[s,t[}).
\]
For general $\gamma\in\mathbb{C}$, the process $(B_{st})$ can be considered as a free analog of the (quantum) Az\'ema martingale with parameter $q=|\gamma|^2$.
\end{example}

\end{document}